\documentclass[a4paper,leqno,12pt]{amsart}
\usepackage{amsmath,amstext,amssymb,amsopn,amsthm,enumitem,float}
\usepackage[utf8]{inputenc}
\usepackage{graphicx}
\usepackage{hyperref}
\usepackage{centernot}
\usepackage{tikz}
\usetikzlibrary{arrows}

\usepackage{pifont}
\newcommand{\cmark}{\ding{51}}%
\newcommand{\xmark}{\ding{55}}%

\newcommand{\NN}{\mathbb{N}}
\newcommand{\RR}{\mathbb{R}}
\newcommand{\MB}{\mathcal{M}^\mathcal{B}}
\newcommand{\BB}{\mathcal{B}}

\theoremstyle{plain}
\newtheorem{theorem}{Theorem}[section]
\newtheorem{lemma}[theorem]{Lemma}
\numberwithin{equation}{section}

\makeatletter
\def\namedlabel#1#2{\begingroup
	#2%
	\def\@currentlabel{#2}%
	\phantomsection\label{#1}\endgroup
}
\makeatother

\begin{document}

\title[$A_\infty$ condition for general bases revisited]{$A_\infty$ condition for general bases revisited:\\ complete classification of definitions}

\author{Dariusz Kosz}
\address{ \newline Dariusz Kosz
	\newline Faculty of Pure and Applied Mathematics
	\newline Wroc\l aw University of Science and Technology 
	\newline Wybrze\.ze Wyspia\'nskiego 27, 50-370 Wroc\l aw, Poland
	\newline \textit{E-mail address:} \textnormal{dariusz.kosz@pwr.edu.pl}		
}


\begin{abstract} We refer to the discussion on different characterizations of the $A_\infty$ class of weights, initiated by Duoandikoetxea, Mart\'in-Reyes, and Ombrosi. Twelve definitions of the $A_\infty$ condition are considered. For cubes in $\RR^d$ every two conditions are known to be equivalent, while for general bases we have a trichotomy: equivalence, one-way implication, or no dependency may occur. In most cases the relations between different conditions have already been established. Here all the unsolved cases are treated and, as a result, a full diagram of the said relations is presented.  
	
	\medskip	
	\noindent \textbf{2010 Mathematics Subject Classification:} Primary 42B25.
	
	\medskip
	\noindent \textbf{Key words:} weights, $A_\infty$ class, maximal operator.
\end{abstract}

\maketitle

\section{Introduction} \label{sec:intro}
While dealing with the theory of weights in the Euclidean setting $\RR^d$, one can find in the literature several statements referred to as the $A_\infty$ condition. The most classical definition is due to Muckenhoupt~\cite{Mu74}. It is said there that a locally integrable function $w \colon \RR^d \to [0,\infty)$ is in the $A_\infty$ class if for each $\epsilon \in (0,1)$ there exists $\delta \in (0,1)$ such that
\[
\Big( |E| < \delta |Q| \Big) \implies \Big( w(E) < \epsilon w(Q) \Big) 
\]
holds, whenever $Q$ is a $d$-dimensional cube and $E$ is its arbitrary measurable subset. This nomenclature has a very natural interpretation as an endpoint for the $A_p$ theory. Indeed, Muckenhoupt showed that $w$ satisfies the above condition if and only if it belongs to the $A_p$ class for some $p \in (1, \infty)$.

At about the same time, Coifman and Fefferman~\cite{CF74} proposed another approach based on verifying the following inequality
\[
\frac{w(E)}{w(Q)} \leq C \bigg( \frac{|E|}{|Q|} \bigg)^\delta,
\]  
where $Q,E$ are as before, while $C, \delta >0$ are universal constants depending only on $w$. It was then proven that the two conditions lead to the same class of weights.

Apart from those mentioned above, there are many other definitions of the $A_\infty$ condition, e.g., the ones using maximal operators (cf.~\cite{Fu78},\cite{Wil87}), medians~(cf.~\cite{OS89}) or a~substitute for the usual $A_p$ condition obtained by letting $p$ to infinity~(cf.~\cite{Hr84},\cite{GR85}). For more detailed information we refer the reader to the survey~\cite{DO85}.

Although the Euclidean theory of weights has been well understood, things got complicated as other settings began to be explored. Indeed, it was observed that some of the equivalences between definitions are consequences of special properties of $\RR^d$ or the geometry of cubes, so that they may no longer be true in other contexts. A natural problem arose --- {\it which definition should be the proper one in general?} --- and the authors have so far not agreed to make any particular version the leading one (compare, for example, \cite{St93},\cite{Duo01},\cite{HP13}).

Recently, Duoandikoetxea, Mart\'in-Reyes, and Ombrosi took a big step towards making the knowledge in this field systematized. In their remarkable paper \cite{DMO16} twelve conditions, which are known to be equivalent in the case of cubes (or balls) in $\RR^d$, are collected and, in the context of general bases associated to arbitrary measure spaces, the exact relations between them are established. In some cases equivalences or one-way implications are proven and in some others suitable counterexamples are constructed. 

Despite its high level of precision, the said article does not cover the problem in full. Indeed, a few cases remained unsolved, as the authors were not able to provide neither a proof nor a counterexample. Thus, it seems very desirable to complete the task, by examining all the missing relations. This is actually what we aim to do here. 

\begin{theorem} \label{thm:1}
	Table~\ref{table} in Section~\ref{sec:main} describes whether an implication of the form \linebreak $\textrm{{\bf (PA)}} \implies \textrm{{\bf (PB)}}$ holds or not for general bases, where $\textrm{{\bf (PA)}}$ and $\textrm{{\bf (PB)}}$ are any two of the twelve $A_\infty$ conditions collected in \cite{DMO16} (for a diagram with the obtained relations see Figure~\ref{figure} in Section~\ref{sec:main}). 
\end{theorem}

The rest of the paper is organized as follows. In Section~\ref{sec:main} we give basic definitions, recall which cases of our problem have already been solved, and present the postulated complete diagram of the studied relations. In Section~\ref{sec:proofs} we prove Theorem~\ref{thm:1}, by treating all the unsolved cases.      

\section{Main result} \label{sec:main}
From now on, we shall deal with a measure space $(X, \Sigma, \mu)$, where $X$ is a set, $\Sigma$ is a $\sigma$-algebra of subsets of $X$, and $\mu \colon \Sigma \to [0,\infty]$ is a $\sigma$-finite measure. A~{\it basis} $\BB$ is a collection of sets $B \in \Sigma$ satisfying $|B| \in (0,\infty)$ (here and later on for $\mu$-measurable sets $E \subset X$ we write simply $|E|$ instead of $\mu(E)$). As in \cite{DMO16}, we assume that the
elements of $\BB$ cover $X$ except for a set of $\mu$-measure zero.

Given a {\it weight} $w$ associated to $\BB$ (that is, a $\mu$-measurable non-negative function such that the integral of $w$ over $B$ is finite for all $B \in \BB$), we introduce the following notation: 
\begin{itemize}
	\item $w(B) := \int_B w \, {\rm d}\mu$ (the {\bf integral} of $w$ over $B$),
	\item $w_B := w(B) / |B|$ (the {\bf average} of $w$ in $B$),
	\item $m(w;B) := \inf \{ t : |x \in B : w(x) > t| < |B| / 2\}$ (the {\bf median} of $w$ in $B$). 
\end{itemize}
Also, we define the {\it maximal operator} associated to $\BB$ by
\[
\MB f(x) := \sup_{x \in B \in \BB} \frac{1}{|B|} \int_B f \, {\rm d}\mu, \qquad x \in X,
\]
where $f$ is any $\mu$-measurable function.

In \cite{DMO16} the following conditions on weights $w$ were mentioned.
\begin{enumerate}
	\item[\namedlabel{P1}{\bf (P1)}] There exists $p \in (1, \infty)$ such that $w$ belongs to the Muckenhoupt class $A_{p, \BB}$, that is, there exists $C_{p, \BB} > 0$ such that
	\[
	\bigg( \frac{1}{|B|} \int_B w \, {\rm d}\mu \bigg) \bigg( \frac{1}{|B|} \int_B w^{1-p'} \, {\rm d}\mu  \bigg)^{p-1} \leq C_{p, \BB}
	\]
	holds for all sets $B \in \BB$, where $\frac{1}{p} + \frac{1}{p'} = 1$. 
	\item[\namedlabel{P1'}{\bf (P1')}] There exist $\delta, C > 0$ such that 
	\[
	\frac{|E|}{|B|} \leq C \bigg( \frac{w(E)}{w(B)} \bigg)^\delta
	\] 
	holds for all sets $B \in \BB$ and all $\mu$-measurable subsets $E \subset B$.
	\item[\namedlabel{P2}{\bf (P2)}] There exists $C > 0$ such that
	\[
	\frac{1}{|B|} \int_{B} w \, {\rm d}\mu \leq C \exp \bigg( \frac{1}{|B|} \int_{B} \log w \, {\rm d}\mu \bigg)
	\]
	holds for all sets $B \in \BB$.
	\item[\namedlabel{P2'}{\bf (P2')}] There exists $C > 0$ such that
	\[
	\frac{1}{|B|} \int_{B} w \, {\rm d}\mu \leq C \bigg( \frac{1}{|B|} \int_{B} w^s \, {\rm d}\mu \bigg)^{1/s}
	\]
	holds for every $s \in (0,1)$ and for all sets $B \in \BB$.
	\item[\namedlabel{P3}{\bf (P3)}] There exists $q \in (1, \infty)$ such that $w$ belongs to the reverse H\"older class $RH_{q, \BB}$, that is, there exists $C_{q, \BB} > 0$ such that
	\[
	 \bigg( \frac{1}{|B|} \int_B w^q \, {\rm d}\mu \bigg)^{1/q} \leq \frac{C_{q, \BB}}{|B|} \int_B w \, {\rm d}\mu
	\]
	holds for all sets $B \in \BB$.
	\item[\namedlabel{P3'}{\bf (P3')}] There exist $\delta, C > 0$ such that 
	\[
	\frac{w(E)}{w(B)} \leq C \bigg( \frac{|E|}{|B|} \bigg)^\delta
	\] 
	holds for all sets $B \in \BB$ and all $\mu$-measurable subsets $E \subset B$.
	\item[\namedlabel{P4}{\bf (P4)}] There exist $\alpha, \beta \in (0,1)$ such that the implication
	\[
	\Big( |E| < \alpha |B| \Big) \implies \Big( \mu(E) \leq \beta \mu(B) \Big)
	\]
	holds for all sets $B \in \BB$ and all $\mu$-measurable subsets $E \subset B$.
	\item[\namedlabel{P4'}{\bf (P4')}] There exist $\alpha, \beta \in (0,1)$ such that
	\[
	\big| \{ x \in B : w(x) \leq \alpha w_B \} \big| \leq \beta |B|
	\]
	holds for all sets $B \in \BB$.
	\item[\namedlabel{P5}{\bf (P5)}] There exists $C > 0$ such that
	\[
	w_B \leq C m(w; B)
	\]
	holds for all sets $B \in \BB$.
	\item[\namedlabel{P6}{\bf (P6)}] There exists $C > 0$ such that
	\[
	\int_B w \log^+ \frac{w}{w_B} \, {\rm d}\mu \leq C w(B)
	\]
	holds for all sets $B \in \BB$.
	\item[\namedlabel{P7}{\bf (P7)}] There exists $C > 0$ such that
	\[
	\int_B \MB (w \chi_B) \, {\rm d}\mu \leq C w(B)
	\]
	holds for all sets $B \in \BB$.
	\item[\namedlabel{P8}{\bf (P8)}] There exist $C, \beta > 0$ such that
	\[
	w\big(\{  x \in B : w(x) \geq \lambda \}\big) \leq C \lambda \, \big|\{  x \in B : w(x) \geq \beta \lambda \}\big|
	\]
	holds for all sets $B \in \BB$ and every $\lambda > w_B$.
\end{enumerate}

We now briefly recall the main results, concerning the relations between these conditions, that were obtained in \cite{DMO16}.  

\smallskip \noindent {\bf Equivalences (\cite[Theorem 3.1]{DMO16}).}
\begin{itemize} 
	\setlength\itemsep{0.1em}
	\item \ref{P1} $\iff$ \ref{P1'}
	\item \ref{P2} $\iff$ \ref{P2'}
	\item \ref{P3} $\iff$ \ref{P3'}
	\item \ref{P4} $\iff$ \ref{P4'}
\end{itemize}

\smallskip \noindent {\bf Implications (\cite[Theorem 4.1]{DMO16}).}
\begin{itemize}
	\setlength\itemsep{0.1em}
	\item \ref{P1} $\implies$ \ref{P2} $\implies$ \ref{P5} $\implies$ \ref{P4}
	\item \ref{P8} $\implies$ \ref{P3} $\implies$ \ref{P6} $\implies$ \ref{P4}
\end{itemize}

\smallskip \noindent {\bf Conditional implications (\cite[Theorem 4.2]{DMO16}).}
\begin{itemize}
	\setlength\itemsep{0.1em}
	\item \ref{P1} $\implies$ \ref{P7} (provided that $\BB$ is a Muckenhoupt basis, that is, the maximal operator $\MB$ is bounded on $L^p(w)$ for each $p \in (1, \infty)$ and for every $w \in A_{p, \BB}$) 
	\item \ref{P3} $\implies$ \ref{P7} (provided that $\MB$ is bounded on $L^p$ for every $p > 1$)
	\item \ref{P6} $\implies$ \ref{P7} (provided that $\MB$ is of weak type $(1,1)$)
	\item \ref{P2} $\implies$ \ref{P7} (provided that $\MB$ is bounded on $L^p$ for some $p > 1$ and that $\MB(w \chi_B)(x) = \sup_{ B \supset B' \in \BB} w_{B'}$ holds for each $x \in B \in \BB$)
\end{itemize}

\smallskip \noindent {\bf Counterexamples (\cite[Counterexamples 1--3 and 5--7]{DMO16}).} All counterexamples are built for $\BB = \{ (0,b) : b > 0\}$ considered as a basis in $(0,\infty)$ equipped with the Lebesgue measure. In particular, all the conditional implications mentioned above are valid.
\begin{itemize}
\setlength\itemsep{0.1em}
\item \ref{P8} $\centernot\implies$ \ref{P5} (thus $\textbf{(PA)} \centernot\implies \textbf{(PB)}$ if $\textbf{A} \in \{{\bf 3,4,6,7,8}\}$ and $\textbf{B} \in \{{\bf 1,2,5}\}$)
\item \ref{P2} $\centernot\implies$ \ref{P1} (thus also \ref{P5} $\centernot\implies$ \ref{P1})
\item \ref{P1} $\centernot\implies$ \ref{P6} (thus $\textbf{(PA)} \centernot\implies \textbf{(PB)}$ if $\textbf{A} \in \{{\bf 1,2,4,5,7}\}$ and $\textbf{B} \in \{{\bf 3,6,8}\}$)
\item \ref{P6} $\centernot\implies$ \ref{P3} (thus also \ref{P6} $\centernot\implies$ \ref{P8})
\item \ref{P3} $\centernot\implies$ \ref{P8}
\item \ref{P7} $\centernot\implies$ \ref{P4}
\end{itemize}
\smallskip

As mentioned before, our aim is to complete the picture, by examining all the missing relations between the above conditions. The final result is presented in Table~\ref{table}. Of course, in view of the four equivalences listed above, we may omit the conditions \ref{P1'}, \ref{P2'}, \ref{P3'}, and \ref{P4'}. 

\begin {table}[H] 
\vspace*{0.3cm}
\caption {Occurrence of implications of the form $\textbf{(PA)}\implies\textbf{(PB)}$ for general bases. The conditions $\textbf{(PA)}$ and $\textbf{(PB)}$ appear in the first column and row, respectively. The symbol $[ \cdot ]$ indicates that the result was not known before.} \label{table}

\begin{center}
	\begin{tabular}{ | c | c | c | c | c | c | c | c | c |}		 
		\hline
		 $\implies$ & \ref{P1} & \ref{P2} & \ref{P3} & \ref{P4} & \ref{P5} & \ref{P6} & \ref{P7} & \ref{P8} \\ \hline
		\ref{P1} & $=$ & \cmark & \xmark & \cmark & \cmark & \xmark & \rm{[}\xmark\rm{]} & \xmark \\ \hline
		\ref{P2} & \xmark & $=$ & \xmark & \cmark & \cmark & \xmark & \rm{[}\xmark\rm{]} & \xmark \\ \hline
		\ref{P3} & \xmark & \xmark & $=$ & \cmark & \xmark & \cmark & \rm{[}\xmark\rm{]} & \xmark \\ \hline
		\ref{P4} & \xmark & \xmark & \xmark & $=$ & \xmark & \xmark & \rm{[}\xmark\rm{]} & \xmark \\ \hline
		\ref{P5} & \xmark & \rm{[}\xmark\rm{]} & \xmark & \cmark & $=$ & \xmark & \rm{[}\xmark\rm{]} & \xmark \\ \hline
		\ref{P6} & \xmark & \xmark & \xmark & \cmark & \xmark & $=$ & \rm{[}\xmark\rm{]} & \xmark \\ \hline
		\ref{P7} & \xmark & \xmark & \xmark & \xmark & \xmark & \xmark & $=$ & \xmark \\ \hline
		\ref{P8} & \xmark & \xmark & \cmark & \cmark & \xmark & \cmark & \rm{[}\xmark\rm{]} & $=$ \\ \hline
	\end{tabular}
\end{center}
\end{table}

In order to prove Theorem~\ref{thm:1} it suffices to show the following three lemmas.

\begin{lemma} \label{lem1}
	There exist a measure space $(X, \Sigma, \mu)$, a basis $\BB \subset \Sigma$, and a weight $w$ associated to $\BB$ such that $w$ satisfies \ref{P5} and does not satisfy \ref{P2}. 
\end{lemma} 

\begin{lemma} \label{lem2}
	There exist a measure space $(X, \Sigma, \mu)$, a basis $\BB \subset \Sigma$, and a weight $w$ associated to $\BB$ such that $w$ satisfies \ref{P1} and does not satisfy \ref{P7}. 
\end{lemma}

\begin{lemma} \label{lem3}
	There exist a measure space $(X, \Sigma, \mu)$, a basis $\BB \subset \Sigma$, and a weight $w$ associated to $\BB$ such that $w$ satisfies \ref{P8} and does not satisfy \ref{P7}. 
\end{lemma}

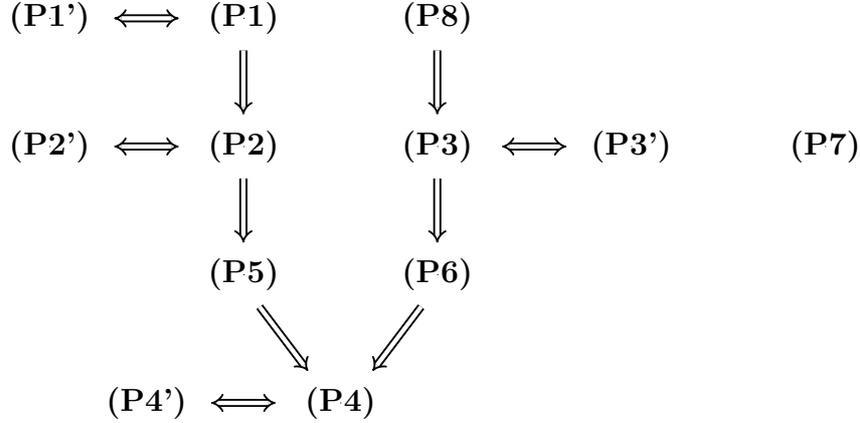
\begin{figure}[H] 
	\begin{tikzpicture}[scale=0.85]
	\draw (4,7) node {\ref{P1}} -- (4,7);
	\draw (4,5) node {\ref{P2}} -- (4,5);
	\draw (7,5) node {\ref{P3}} -- (7,5);
	\draw (5.5,1) node {\ref{P4}} -- (5.5,1);
	\draw (4,3) node {\ref{P5}} -- (4,3);
	\draw (7,3) node {\ref{P6}} -- (7,3);
	\draw (13,5) node {\ref{P7}} -- (13,5);
	\draw (7,7) node {\ref{P8}} -- (7,7);
	\draw (1,7) node {\ref{P1'}} -- (1,7);
	\draw (1,5) node {\ref{P2'}} -- (1,5);
	\draw (10,5) node {\ref{P3'}} -- (10,5);
	\draw (2.5,1) node {\ref{P4'}} -- (2.5,1);
	\draw[implies-implies,double equal sign distance, line width=0.25mm] (2,7) -- (3,7);
	\draw[implies-implies,double equal sign distance, line width=0.25mm] (2,5) -- (3,5);
	\draw[implies-implies,double equal sign distance, line width=0.25mm] (8,5) -- (9,5);
	\draw[implies-implies,double equal sign distance, line width=0.25mm] (3.5,1) -- (4.5,1);
	\draw[-implies,double equal sign distance, line width=0.25mm] (4,6.5) -- (4,5.5);
	\draw[-implies,double equal sign distance, line width=0.25mm] (4,4.5) -- (4,3.5);
	\draw[-implies,double equal sign distance, line width=0.25mm] (4.25,2.5) -- (5,1.5);
	\draw[-implies,double equal sign distance, line width=0.25mm] (7,6.5) -- (7,5.5);
	\draw[-implies,double equal sign distance, line width=0.25mm] (7,4.5) -- (7,3.5);
	\draw[-implies,double equal sign distance, line width=0.25mm] (6.75,2.5) -- (6,1.5);
	\draw (0,0) node {} -- (0,0);
	\draw (0,7.5) node {} -- (0,7.5);
	\end{tikzpicture}
	\caption{Relations between different $A_\infty$ conditions in the context of general bases.}
	\label{figure}
\end{figure}

For the reader's convenience, in Figure~\ref{figure} we present the final result as a diagram of relations between different $A_\infty$ conditions. 

\section{Proofs of Lemmas~\ref{lem1}--\ref{lem3}} \label{sec:proofs}
In this section we present the proofs of Lemmas~\ref{lem1}--\ref{lem3}. We shall work with countable discrete spaces $X$ in order to provide examples which are as simple as possible (in particular, $\Sigma$ is always the family of all subsets of $X$). However, we emphasize that a suitable continuous counterpart can be assigned to each of these examples. Indeed, given $(X, \Sigma, \mu)$, $\BB$, and $w$, the discrete objects specified in one of the proofs, one define the product space
\[
(\overline{X}, \overline{\Sigma}, \overline{\mu} ) := \big(X \times [0,1), \Sigma \times \BB_{[0,1)}, \mu \times \lambda_{[0,1)} \big)
\]   
(here the symbols $\BB_{[0,1)}$ and $\lambda_{[0,1)}$ stay for the $\sigma$-algebra of Borel subsets of $[0,1)$ and the Lebesgue measure on $[0,1)$, respectively), the basis
\[
\overline{\BB} := \{ B \times [0,1) : B \in \BB \},
\]
and the weight
\[
\overline{w}(x,t) := w(x), \qquad x \in X, \ t \in [0,1).
\]
Then it turns out that all the properties we want $w$ to satisfy are also satisfied by $\overline{w}$. 

In fact, even more can be said. Assuming $X = \{ x_1, x_2, \dots \}$ and $|\{ x_j \}| = \alpha_j > 0$, $j \in \mathbb{N}$, with $\sum_{j \in \mathbb{N}} \alpha_j = \infty$ (this will indeed be satisfied in each of the three examples), one can use the bijective map $\tau \colon \overline{X} \to [0,\infty)$ given by
\[
\tau\big((x_j,t)\big) := \sum_{i < j} \alpha_i + t \alpha_j,
\]
in order to obtain the desired example with the underlying space being just the half-line $[0,\infty)$ equipped with the Lebesgue measure $\lambda_{[0,\infty)}$. In this situation, the only peculiarities are the strange looking basis $ \BB_{[0,\infty)} := \{ \tau^{-1}(B) : B \in \overline{\BB} \}$ and weight $w_{[0,\infty)}(s) := \overline{w} ( \tau^{-1} (s))$. This observation reveals that the structure of the basis plays the main role in the problem, not the specific properties of the measure space such as the doubling condition, for instance.    

In what follows we write $|  x |$ instead of $| \{ x \} |$ (thus, e.g., $w(B) = \sum_{x \in B} w(x) \cdot |x|$). 

\smallskip \noindent
{\bf Proof of Lemma~\ref{lem1}.}
Set
\[
X := \big\{ x_{n,i} : n \in \NN, \ i \in \{0, 1 \} \big\}
\]
and define $\mu$ by letting
\[
|  x_{n,i} | := \left\{ \begin{array}{rl}
2 & \textrm{if } i=0, \\
1 & \textrm{otherwise.} \end{array} \right.
\]
Consider the basis
\[
\BB := \{ B_n : n \in \NN \}, 
\]
where $B_n := \{x_{n,0}, x_{n,1} \}$. We shall show that $w$ defined by
\[
w(x_{n,i}) := \left\{ \begin{array}{rl}
n & \textrm{if } i=0, \\
1 & \textrm{otherwise,} \end{array} \right.
\]
satisfies \ref{P5} and does not satisfy \ref{P2}. 

We first show that \ref{P5} is satisfied with $C = 1$. Indeed, for $n \in \NN$ we have
\[
w_{B_n} = \frac{2n + 1}{2 + 1} \leq n = m(w; B_n). 
\]
In order to show that \ref{P2} cannot hold observe that
\begin{align*}
\frac{\sum_{x \in B_n} w(x) \cdot |x|}{|B_n|}  = \frac{2n+1}{3} & > \frac{2n^{1/3}}{3} \exp \bigg( \frac{2 \log n + \log 1}{3} \bigg) \\ 
& = \frac{2n^{1/3}}{3} \exp \bigg( \frac{\sum_{x \in B_n} \log w(x) \cdot |x|}{|B_n|} \bigg).
\end{align*}

\smallskip \noindent
{\bf Proof of Lemma~\ref{lem2}.} 
Set
\[
X := \big\{ x_{n,i} : n \in \NN, \ i \in \{0, 1, 2, \dots, 4^n \} \big\}
\]
and define $\mu$ by letting
\[
|  x_{n,i} | := \left\{ \begin{array}{rl}
1 & \textrm{if } i=0, \\
2^{-n} & \textrm{otherwise.} \end{array} \right.
\]
Consider the basis
\[
\BB := \big\{ B_{n,i} : n \in \NN, \ i \in \{0, 1, 2, \dots, 4^n \} \big\}, 
\]
where
\[
B_{n,i} := \left\{ \begin{array}{rl}
\{ x_{n,0}, x_{n,1}, \dots, x_{n,4^n} \} & \textrm{if } i=0, \\
\{ x_{n,0}, x_{n,i} \} & \textrm{otherwise.} \end{array} \right.
\]
We shall show that $w$ defined by
\[
w(x_{n,i}) := \left\{ \begin{array}{rl}
1 & \textrm{if } i=0, \\
2^{-n} & \textrm{otherwise,} \end{array} \right.
\]
satisfies \ref{P1} and does not satisfy \ref{P7}.

We first show that \ref{P1} holds, by proving that \ref{P1'} is satisfied with $C=4$ and $\delta = \frac{1}{2}$. Take $B_{n,i} \in \BB$ and non-empty $E \subset B_{n,i}$. We consider two cases, $i=0$ and $i > 0$. If $i=0$, then $|B| = 1+2^n$ and $w(B) = 2$. Moreover, we have $|E| = j + k 2^{-n}$ and $w(E) = j + k 4^{-n}$ for some $j \in \{0,1\}$ and $k \in \{0,1, \dots, 4^n\}$ such that $j+k > 0$. Thus,
\[
\frac{|E|}{|B|} < \frac{j}{2^n} + \frac{k}{4^n} \leq \frac{w(E)}{w(B)} + 2 \cdot \frac{w(E)}{w(B)} \leq 3 \cdot \bigg( \frac{w(E)}{w(B)} \bigg)^{1/2}.
\]
On the other hand, if $i>0$, then $|B| = 1 + 2^{-n}$ and $w(B) = 1+4^{-n}$. Moreover, we have $|E| = j + k 2^{-n}$ and $w(E) = j + k 4^{-n}$ for some $j,k \in \{0,1\}$ such that $j+k > 0$. Thus,
\[
\frac{|E|}{|B|} < j + \frac{k}{2^n} \leq 2 \cdot \frac{w(E)}{w(B)} + 2 \cdot \bigg( \frac{w(E)}{w(B)} \bigg)^{1/2} \leq 4 \cdot \bigg( \frac{w(E)}{w(B)} \bigg)^{1/2}.
\]
Consequently, \ref{P1} is satisfied.

In order to show the second part observe that
\[
\MB (w \chi_{B_{n,0}})(x_{n,i}) \geq \frac{w(B_{n,i})}{|B_{n,i}|} = \frac{1+4^{-n}}{1+2^{-n}} > \frac{1}{2}, \qquad i \in \{1,2, \dots, 4^n\}.
\]
Thus,
 \[
 \sum_{x \in B_{n,0}} \MB (w \chi_{B_{n,0}})(x) \cdot |x| > 4^n \cdot \frac{1}{2} \cdot \frac{1}{2^n} = 2^{n-2} w(B_{n,0})
 \]
and, consequently, \ref{P7} cannot hold.

\smallskip \noindent
{\bf Proof of Lemma~\ref{lem3}.}
Let $m_0 = 0$ and $m_n = 4 + 4^2 + \dots + 4^n$ for $n \in \NN$. Set
\[
X := \big\{ x_{n,i} : n \in \NN, \ i \in \{0, 1, 2, \dots, m_{n} \} \big\},
\]
and specify $\mu$ to be the counting measure. Consider the basis
\[
\BB := \big\{ B_{n,i} : n \in \NN, \ i \in \{0, 1, 2, \dots, m_n \} \big\}, 
\]
where
\[
B_{n,i} := \left\{ \begin{array}{rl}
\{ x_{n,0}, x_{n,1}, \dots, x_{n,m_n} \} & \textrm{if } i=0, \\
\{ x_{n,0}, x_{n,i} \} & \textrm{otherwise.} \end{array} \right.
\]
We shall show that $w$ defined by
\[
 w(x_{n,i}) := \left\{ \begin{array}{rl}
1 & \textrm{if } i=0, \\
2^{-j} & \textrm{if } m_{j-1} < i \leq m_j \textrm{ for } j \in \{1, \dots, n\}. \end{array} \right.
\]
satisfies \ref{P8} and does not satisfy \ref{P7}. 

We first show that \ref{P8} is satisfied with $C=4$ and $\beta = 1$. Take $B_{n,i} \in \BB$. We consider two cases, $i=0$ and $i > 0$. If $i=0$, then $w_{B_{n,i}} > 2^{-n}$. Thus, if $2^{-k-1} < \lambda \leq 2^{-k}$ holds for some $k \in \{0, \dots, n-1\}$, then
\[
w\big(\{  x \in B_{n,i} : w(x) \geq \lambda \}\big) = 2^0 + 2^1 + \dots + 2^k < 2 \cdot 2^k
\]
and
\[
\big|\{  x \in B_{n,i} : w(x) \geq \lambda \}\big| = 4^0 + 4^1 + \dots + 4^k > 4^k.
\]
Combining the above, we get
\[
w\big(\{  x \in B_{n,i} : w(x) \geq \lambda \}\big) < 2 \cdot 2^k = 4 \cdot 2^{-k-1} \cdot 4^k < 4 \lambda \, \big|\{  x \in B_{n,i} : w(x) \geq \lambda \}\big|.
\]
On the other hand, if $i>0$, then $w_{B_{n,i}} > \frac{1}{2}$. Thus, for $\frac{1}{2} < \lambda \leq 1$ we have
\[
w\big(\{  x \in B_{n,i} : w(x) \geq \lambda \}\big) = 1 < 2 \lambda = 2 \lambda \, \big|\{  x \in B_{n,i} : w(x) \geq \lambda \}\big|.
\]
Consequently, \ref{P8} is satisfied.

In order to show the second part observe that
\[
\MB (w \chi_{B_{n,0}})(x_{n,i}) \geq \frac{w(B_{n,i})}{|B_{n,i}|} > \frac{w(x_{n,0}) \cdot |x_{n,0}|}{2} = \frac{1}{2}, \qquad i \in \{1,2, \dots, m_n\}.
\]
Thus,
\[
\sum_{x \in B_{n,0}} \MB (w \chi_{B_{n,0}})(x) \cdot |x| > \frac{m_n}{2} > \frac{4^n}{2} = \frac{2^n}{4} \cdot 2 \cdot 2^n > \frac{2^n}{4} w(B_{n,0})
\]
and, consequently, \ref{P7} cannot hold.

\end{document}